\documentclass[10pt,twocolumn]{IEEEtran}

\usepackage{graphics} 
\usepackage{epsfig} 
\usepackage{amsmath} 
\usepackage{amssymb}  

\usepackage{graphicx,xcolor}
\usepackage{stfloats}
\usepackage{float}

\graphicspath{{./Figs/}}

\providecommand{\tabularnewline}{\\}
\usepackage{caption}

\newcommand*{\TitleFont}{%
      \fontsize{16}{20}%
      \selectfont}

\begin{document}

\title{\vspace{-12pt}\TitleFont Assessing the Economics of Customer-Sited Multi-Use Energy Storage \vspace{-3pt}}

\author{Wuhua Hu, Ping Wang, and Hoay Beng Gooi
\thanks{This work was supported in part by the Energy Innovation Programme Office (EIPO) through the National Research Foundation and Singapore Economic Development Board. W. Hu is with the Signal Processing Department, Institute for Infocomm Research, A*STAR, Singapore. P. Wang is with the School of Computer Engineering and H. B. Gooi is with the School of Electrical and Electronic Engineering, Nanyang Technological University, Singapore. E-mails: {\tt\small huwh@i2r.a-star.edu.sg, \{wangping, ehbgooi\}@ntu.edu.sg}.}
\vspace{-12pt}}

\maketitle
\thispagestyle{empty}
\pagestyle{plain}

\begin{abstract}
This paper presents an approach to assess the economics of
customer-sited energy storage systems (ESSs) which are owned and operated by a customer.
The ESSs can participate in frequency regulation and spinning reserve
markets, and are used to help the customer consume available renewable energy and reduce electricity bill. A rolling-horizon approach is developed to optimize the service schedule, and the resulting costs and revenues are used to assess economics of the ESSs. The economic assessment approach is illustrated with case studies,
from which we obtain some new observations on profitability of
the customer-sited multi-use ESSs.
\end{abstract}


\section{Introduction}

Energy storage systems (ESSs) are a promising ingredient
for reliable integration of renewable energies into future
power grids. ESSs are however costly at the present stage, and recent studies
showed that they are unlikely to generate a net profit if ESSs are
used to provide a single service. This motivates the use of ESSs for
multiple service provision \cite{Fitzgerald2015}.

When ESSs are scheduled for multiple services concurrently,
potential conflicts occur due to the limited power and energy capacities
available. An ideal scheduling approach needs to address the conflicts
in an optimal way, such that the net profit is maximized subject to operational and service constraints. So far, only a few studies have
been conducted, partially addressing the encountered challenges.
Among them, \cite{you2010economic} presents a coarse framework to
investigate the net profit, but leaves the nontrivial modeling of the optimization objective and
constraints to readers for specific applications. Reference \cite{moreno2015milp}
develops a mixed-integer linear programming model focusing
on maximizing the revenue without
considering the costs of ESSs. More recently, \cite{wu2015energy}
presents a concrete optimization formulation in a rolling-horizon framework. However, the formulation does not appropriately capture the operating costs of ESSs which are dependent on
their varying charge and discharge rates \cite{trippe2014charging}.

This work considers customer-sited ESSs which provide multiple services.
The ESSs are used to participate in
regulation and spinning reserve markets, and help the customer consume available renewable energy and reduce time-of-use (TOU) electricity bill. We develop a comprehensive
scheduling model which captures the dynamics of ESSs and associated aging costs, the supported services and associated revenues, and all major service and operational constraints. By optimizing the schedule using a rolling-horizon approach, we are able to assess profitability of the ESSs. Different from aforementioned literature \cite{you2010economic,moreno2015milp,wu2015energy},
we include the support of self-consumption of renewable
energy and embed a more realistic aging model for the ESSs. When the storage is made of Li-ion batteries, the aging
model characterizes the battery aging cost in terms of its instant
charge/discharge rate and the duration, which was experimentally established
in \cite{trippe2014charging}.

\section{Modeling the ESSs and Their Services}  \label{sec: Modeling-the-ESSs}

We use $\mathcal{N}$ to denote the set of ESSs.
The time is discretized into slots, each with a duration of $T_{{\rm s}}$. The charge and discharge of ESSs are scheduled periodically to support self-consumption of renewable energy, frequency regulation, spinning reserve, and TOU electric bill reduction. The mathematical
models of the ESSs and the four supported services are developed in
this section.

\subsection{Modeling the ESSs}

We assume that a customer owes and operates multiple ESSs, each of which
follows a generic model used in \cite{hu2016towards}. Let the charge
and discharge rates of ESS $i$ be scheduled as $p_{i,t}^{{\rm c}}$
and $p_{i,t}^{{\rm d}}$ for time slot $t$, respectively. And let
$v_{i,t}^{{\rm c}}$ indicate the working mode of the ESS $i$, which
is 1 (or 0) if it is \textit{not} discharged (or \textit{not} charged).
These variables satisfy
\begin{equation}
0\le p_{i,t}^{{\rm c}}\le v_{i,t}^{{\rm c}}p_{i,\max}^{{\rm c}},\thinspace\thinspace0\le p_{i,t}^{{\rm d}}\le(1-v_{i,t}^{{\rm c}})p_{i,\max}^{{\rm d}},\label{cons: ES charge-discharge limits}
\end{equation}
where $p_{i,\max}^{{\rm c}}$ and $p_{i,\max}^{{\rm d}}$ are the
corresponding upper bounds. The two constraints ensure that charge
and discharge comply with the rate limits and do not happen in the same time slot.
After charge/discharge, the state of charge (SOC) of the ESS $i$,
denoted by $s_{i,t}$, renews into
\begin{equation}
s_{i,t}=s_{i,t-1}+T_{s}(\eta_{i}^{{\rm c}}p_{i,t}^{{\rm c}}-p_{i,t}^{{\rm d}}/\eta_{i}^{{\rm d}})/E_{i}^{\text{cap}},\label{cons: SOC update}
\end{equation}
where $\eta_{i}^{{\rm c}},\eta_{i}^{{\rm d}}\in(0,1)$ are the energy
conversion coefficients, and $E_{i}^{{\rm cap}}$ is the energy capacity
of ESS $i$. The SOC must be maintained within certain limits in order to
protect the ESSs, and this will be discussed later in Section \ref{sub: Service-support-constraint}.

Both charge and discharge incur an aging cost, which is the money loss
of the initial investment. Let the cost be estimated as $C_{i}^{{\rm c}}(p_{i,t}^{{\rm c}})$
and $C_{i}^{{\rm d}}(p_{i,t}^{{\rm d}})$ for charging and discharging ESS $i$ at the rates of $p_{i,t}^{{\rm c}}$ and $p_{i,t}^{{\rm d}}$
for one hour, respectively. The cost of operating ESS $i$ in
time slot $t$ is then given by
\[
C_{i}(p_{i,t}^{{\rm c}},p_{i,t}^{{\rm d}})=T_{{\rm s}}C_{i}^{{\rm c}}(p_{i,t}^{{\rm c}})+T_{{\rm s}}C_{i}^{{\rm d}}(p_{i,t}^{{\rm d}}).
\]

If the ESSs use Li-ion batteries, the cost can be approximated by
a piece-wise linear function which is further obtained by solving
the following linear program \cite{trippe2014charging,hu2016towards}:{
\par \vspace{-6pt}
 {\small{}
\begin{equation}
\begin{aligned} & C_{i}(p_{i,t}^{{\rm c}},p_{i,t}^{{\rm d}})\approx\frac{\alpha_{i}T_{{\rm s}}}{0.8E_{i}^{\text{cap}}}\min_{\zeta_{i,t}^{{\rm ESS}}}\zeta_{i,t}^{{\rm ESS}},\\
&\text{s.t., }\gamma_{i}\eta_{i}^{{\rm c}}[1000\times a_{k}^{{\rm ESS}}(p_{i,t}^{{\rm c}})^{2}+n_{i}b_{k}^{{\rm ESS}}p_{i,t}^{{\rm c}}]+\frac{1-\gamma_{i}}{\eta_{i}^{{\rm d}}}\\
& \times [1000\times a_{k}^{{\rm ESS}}(p_{i,t}^{{\rm d}})^{2}+n_{i}b_{k}^{{\rm ESS}}p_{i,t}^{{\rm d}}] \le\zeta_{i,t}^{{\rm ESS}},\thinspace\forall\, k\in\mathcal{K}_{{\rm ESS}},
\end{aligned}
\label{cons: C_ESS}
\end{equation}
}}where $\alpha_{i}$ is the unit capital cost (\$/Wh) to purchase ESS
$i$; $\zeta_{i,t}^{\rm ESS}$ is an auxiliary variable; $\gamma_{i}$
is the fraction of a single cyclic aging cost incurred by fully charging
the battery from empty; $n_{i}\triangleq E_{i}^{\text{cap}}/0.0081$,
which is the number of battery modules that form the ESS $i$, each
with a capacity of 0.0081 kW; and $\{a_{k}^{{\rm ESS}},b_{k}^{{\rm ESS}}\}_{k\in\mathcal{K}_{{\rm ESS}}}$
are the coefficients associated with the linear segments as indicated
by a certain set $\mathcal{K}_{{\rm ESS}}$.

\subsection{Service for self-consumption of renewable energy}

The customer has installed renewable energy generators. The aggregate generation power for time slot $t$ is denoted as $p_{t}^{{\rm re}}$. For time slot $t$, let the customer be scheduled to consume the renewable energy at a
rate of $p_{t}^{{\rm re,sc}}$, and the surplus
renewable energy be charged to ESS $i$ at a rate of $p_{i,t}^{{\rm re,c}}$,
and the remaining renewable energy be exported to the market at a
rate of $p_{t}^{{\rm re,s}}$. These power variables satisfy
{\par \vspace{-8pt}
\begin{equation} {\small
\begin{gathered}0\le p_{t}^{{\rm re,sc}}\le d_{t},\,\,0\le p_{i,t}^{{\rm re,c}}\le p_{i,\max}^{{\rm c}},\,\,0\le p_{t}^{{\rm re,s}}\le p_{\max}^{{\rm s}},\\
p_{t}^{{\rm re,sc}}+p_{t}^{{\rm re,s}}+\sum_{i\in\mathcal{N}}p_{i,t}^{{\rm re,c}}\le p_{t}^{{\rm re}},
\end{gathered} }
\label{cons: sc of renewable energy}
\end{equation}
}%
where $d_{t}$ is the load demand of the customer, and
$p_{\max}^{{\rm s}}$ is the maximum power that can be injected to the grid.
The last inequality admits curtailment of surplus renewable generation,
if any.

Given the electricity purchase price $c_{t}^{{\rm p}}$ and sale price $c_{t}^{{\rm s}}$, we can compute
the revenue of consuming renewable energy with the
help of ESSs as{

\par \vspace{-8pt}
{\footnotesize{}
\begin{align*}
R_{{\rm sc}}(p_{t}^{{\rm re,sc}},p_{t}^{{\rm re,s}},\{p_{i,t}^{{\rm re,c}}\}_{i\in\mathcal{N}}) & =T_{{\rm s}}c_{t}^{{\rm p}}(p_{t}^{{\rm re,sc}}+\sum_{i\in\mathcal{N}}p_{i,t}^{{\rm re,c}})+T_{{\rm s}}c_{t}^{{\rm s}}p_{t}^{{\rm re,s}},
\end{align*}
}}%
of which the first part owes to the avoided purchase of energy
from the market, and the second part owes to the surplus renewable
energy exported to the market.

\subsection{Service for frequency regulation}

Frequency regulation aims at stabilizing the grid frequency at a desired
value. Let $u_{t}^{{\rm fr},{\rm up}}$ be an indicator
which is 1 for ramp up regulation and 0 otherwise. The minimum power
to participate in the regulation market is required to be $p_{\min}^{{\rm fr}}$.
The ESSs may participate in the market or
not, as indicated by $v_{t}^{{\rm fr}}$ equal to 1 and 0, respectively.
Let ESS $i$ charge at a rate of $p_{i,t}^{{\rm fr,c}}$
if $u_{t}^{{\rm fr},{\rm up}}=0$ and discharge at a rate of $p_{i,t}^{{\rm fr,d}}$
if $u_{t}^{{\rm fr},{\rm up}}=1$. The power variables satisfy
\begin{equation}
\begin{gathered}0\le p_{i,t}^{{\rm fr,d}}\le v_{t}^{{\rm fr}} u_{t}^{{\rm fr},{\rm up}}p_{i,\max}^{{\rm d}},\thinspace\thinspace\forall\thinspace i\in\mathcal{N},\\
0\le p_{i,t}^{{\rm fr,c}}\le v_{t}^{{\rm fr}} (1-u_{t}^{{\rm fr},{\rm up}})p_{i,\max}^{{\rm c}}\thinspace\thinspace\forall\thinspace i\in\mathcal{N},\\
p_{t}^{{\rm fr}}\triangleq\sum_{i\in\mathcal{N}}(1-u_{t}^{{\rm fr},{\rm up}})p_{i,t}^{{\rm fr,c}}+u_{t}^{{\rm fr},{\rm up}}p_{i,t}^{{\rm fr,d}}\ge v_{t}^{{\rm fr}}p_{\min}^{{\rm fr}},
\end{gathered}
\label{cons: frequency regulation}
\end{equation}
where the first two inequalities ensure that charge and discharge for the regulation service do not happen concurrently.

Consider the payment scheme implemented by PJM in
USA \cite{avendano2015financial,manual2015energy}. The regulation service is paid by the committed
power capacity (denoted by $p_{t}^{{\rm fr}}$) and the regulation
performance (dictated by the performance score $\rho_{t}^{{\rm fr}}$
and the regulation mileage ratio $\mu_{t}^{{\rm fr}}$). The performance
score ($\rho_{t}^{{\rm fr}}$) is computed based on the regulation performance in the past period; and the
mileage ratio ($\mu_{t}^{{\rm fr}}$) is the mileage of the fast regulation
signal divided by the mileage of the slow (or conventional) regulation
signal, both in the past service period. The Regulation Market Capacity Clearing
Price (RMCCP) is denoted by $c_{t}^{{\rm RMCCP}}$, and the Regulation
Market Performance Clearing Price (RMPCP) is denoted by $c_{t}^{{\rm RMPCP}}$. Both prices
are updated at a period of $T_{{\rm s}}$.

The revenue of the regulation service is then computed
as{
\par \vspace{-6pt}
{\footnotesize{}
\begin{align*}
R_{{\rm fr}}(p_{t}^{{\rm fr}}) & =T_{{\rm s}}\rho_{t}^{{\rm fr}}p_{t}^{{\rm fr}}(c_{t}^{{\rm RMCCP}}+c_{t}^{{\rm RMPCP}}\mu_{t}^{{\rm fr}})+T_{{\rm s}}c_{t}^{{\rm p}}\sum_{i\in\mathcal{N}}(p_{i,t}^{{\rm fr,d}}-p_{i,t}^{{\rm fr,c}}),
\end{align*}
}}%
where the first term owes to the service provided, and the second
term accounts for the revenue obtained from the energy charged/discharged to/from
the ESSs.

\subsection{Service as spinning reserve}

Consider a spinning reserve market which periodically publishes a
reserve availability price, denoted by $c_{t}^{{\rm sr}}$. The minimum participation power is required to be $p_{\min}^{{\rm sr}}$,
and the minimum commission time is $T_{\min}^{{\rm sr}}$. The ESSs may be scheduled to support this service, which
is dictated by a binary variable $v_{t}^{{\rm sr}}$, 1 for participation
and 0 otherwise. Let $p_{i,t}^{{\rm sr},{\rm d}}$ be the power reserved by ESS $i$, which
is the commissioned maximum discharge rate under contingencies. The reserved
power satisfies
\begin{equation}
0\le p_{i,t}^{{\rm sr,d}}\le v_{t}^{{\rm sr}}p_{i,\max}^{{\rm d}},\quad\sum_{i\in\mathcal{N}}p_{i,t}^{{\rm sr,d}}\ge v_{t}^{{\rm sr}}p_{\min}^{{\rm sr}}.\label{cons: operating reserve}
\end{equation}
The minimum support time will be enforced via
constraint (\ref{cons: service linkage c}) ahead. The revenue
for this service is calculated by
\[
R_{{\rm sr}}(\{p_{i,t}^{{\rm sr},{\rm d}}\}_{i\in\mathcal{N}})=c_{t}^{{\rm sr}}T_{{\rm s}}\sum_{i\in\mathcal{N}}p_{i,t}^{{\rm sr},{\rm d}}.
\]

\subsection{Service for TOU electricity bill reduction and preparation for future services}

With TOU electricity pricing information, the ESSs
may be used to reduce the electricity bill by charging and discharging the storage appropriately. At the meantime, the ESSs may be charged/discharged
to prepare for future services. Let
the aggregate charge and discharge rates of ESS $i$ for such purposes be scheduled
as $p_{i,t}^{{\rm fs,c}}$ and $p_{i,t}^{{\rm br+fs,d}}$ for time slot $t$, respectively. The revenue obtained from
the charged and discharged energy is computed as{
\par  \vspace{-10pt}{\small{}
\[
R_{{\rm br}}(\{p_{i,t}^{{\rm br+fs,d}},p_{i,t}^{{\rm fs,c}}\}_{i\in\mathcal{N}})=T_{{\rm s}}c_{t}^{{\rm p}}\sum_{i\in\mathcal{N}}(p_{i,t}^{{\rm br+fs,d}}-p_{i,t}^{{\rm fs,c}}),
\]
}}%
which owes to the avoided or desired purchase of energy from the
market. The charge and discharge rates satisfy
\begin{equation}
0\le p_{i,t}^{{\rm br+fs,d}}\le p_{i,\max}^{{\rm d}},\quad0\le p_{i,t}^{{\rm fs,c}}\le p_{i,\max}^{{\rm c}}.\label{cons: bill reduction}
\end{equation}

\subsection{Feasibility constraints to support the multiple services \label{sub: Service-support-constraint}}

To support the four services above, we must have
\begin{gather}
p_{i,t}^{{\rm re,c}}+p_{i,t}^{{\rm fr,c}}+p_{i,t}^{{\rm fs,c}}=p_{i,t}^{{\rm c}}\le v_{i,t}^{{\rm c}}p_{i,\max}^{{\rm c}},\label{cons: service linkage a}\\
p_{i,t}^{{\rm fr,d}}+p_{i,t}^{{\rm br+fs,d}}=p_{i,t}^{{\rm d}}\le(1-v_{i,t}^{{\rm c}})(p_{i,\max}^{{\rm d}}-p_{i,t}^{{\rm sr,d}}),\label{cons: service linkage b}\\
s_{i,\min}+p_{i,t}^{{\rm sr,d}}T_{\min}^{{\rm sr}}/E_{i}^{{\rm cap}}\le s_{i,t}\le s_{i,\max},\label{cons: service linkage c}
\end{gather}
for each $i\in\mathcal{N}$. Constraints (\ref{cons: service linkage a}) and (\ref{cons: service linkage b}) are related to the aggregate charge rate and discharge rate for multiple services, respectively. Constraint (\ref{cons: service linkage c}) imposes SOC limits to protect the ESSs from being over charged
or discharged, in which $s_{i,\min},s_{i,\max}\in(0,1)$ are the required
limits and $p_{i,t}^{{\rm sr,d}}T_{\min}^{{\rm sr}}$ is the energy committed as spinning reserve. The three constraints link
up the four services provided by the ESSs, through which
the conflicts in between will be resolved via optimization.

The right hand side of the inequality in (\ref{cons: service linkage b})
contains a term $v_{i,t}^{{\rm c}}p_{i,t}^{{\rm sr,d}}$, which is
bilinear in the two decision variables. It is desirable to reformulate
this term into a linear equivalent form. Introduce an auxiliary
variable $z_{i,t}$. Then, $z_{i,t}$ is equal to $v_{i,t}^{{\rm c}}p_{i,t}^{{\rm sr,d}}$
if it satisfies the following constraints:
\begin{gather}
\begin{gathered}0\le z_{i,t}\le p_{i,\max}^{{\rm d}}v_{i,t}^{{\rm c}},\\
p_{i,t}^{{\rm sr,d}}+p_{i,\max}^{{\rm d}}(v_{i,t}^{{\rm c}}-1)\le z_{i,t}\le p_{i,t}^{{\rm sr,d}},
\end{gathered}
\label{cons: auxiliary z}
\end{gather}
The equivalence is easy to verify with the McCormick linearization
method \cite{mccormick1976computability}. Therefore, we can replace $v_{i,t}^{{\rm c}}p_{i,t}^{{\rm sr,d}}$
with $z_{i,t}$ subjected to the above constraints, which
makes constraint (\ref{cons: service linkage b}) completely linear
in the variables.

\section{The Storage Management Optimization Problem and Its Solution} \label{sec: The-Storage-Management-problem}

The storage management optimization problem is defined to maximize
the total net profit (namely, minimize the total net loss) over a rolling
horizon subject to service requirements and operational constraints. Given a decision time point $t$, the time slots
within a look-ahead horizon for a size of $H$ are denoted by the set $\mathcal{H}_{t}$. The optimization will be performed using the forecast data over the horizon.

The total net profit (TNP) over the horizon sums up the revenues earned
from the four services minus the operating cost incurred
to the ESSs. With the storage aging cost $C_{i}(p_{i,t}^{{\rm c}},p_{i,t}^{{\rm d}})$
computed from the linear program (\ref{cons: C_ESS}), the TNP can be shown to have
the following specific form:{
\par \vspace{-6pt}
{\scriptsize{}
\begin{align*}
 & \text{TNP}(t)=\sum_{\tau\in\mathcal{H}_{t}}\left(\begin{array}{c}
R_{{\rm sc}}(p_{\tau}^{{\rm re,sc}},p_{t}^{{\rm re,s}},\{p_{i,\tau}^{{\rm re,c}}\}_{i\in\mathcal{N}})\\
+R_{{\rm fr}}(\{p_{i,\tau}^{{\rm fr,c}},p_{i,\tau}^{{\rm fr,d}}\}_{i\in\mathcal{N}})+R_{{\rm sr}}(\{p_{i,\tau}^{{\rm sr},{\rm d}}\}_{i\in\mathcal{N}})\\
+R_{{\rm br}}(\{p_{i,\tau}^{{\rm br+fs,d}},p_{i,\tau}^{{\rm fs,c}}\}_{i\in\mathcal{N}})\\
-\sum_{i\in\mathcal{N}}C_{i}(p_{i,\tau}^{{\rm c}},p_{i,\tau}^{{\rm d}})
\end{array}\right)\\
 & =T_{{\rm s}}\sum_{\tau\in\mathcal{H}_{t}}\left(\begin{array}{c}
c_{\tau}^{{\rm p}}p_{\tau}^{{\rm re,sc}}+c_{\tau}^{{\rm s}}p_{\tau}^{{\rm re,s}}+\rho_{\tau}^{{\rm fr}}p_{\tau}^{{\rm fr}}(c_{\tau}^{{\rm RMCCP}}+c_{\tau}^{{\rm RMPCP}}\mu_{\tau}^{{\rm fr}})\\
+\sum_{i\in\mathcal{N}}\left(\begin{array}{c}
c_{\tau}^{{\rm p}}(p_{i,\tau}^{{\rm d}}-p_{i,\tau}^{{\rm c}}+2p_{i,\tau}^{{\rm re,c}})\\
+c_{\tau}^{{\rm sr}}p_{i,\tau}^{{\rm sr},{\rm d}}-\min_{\zeta_{i,\tau}^{{\rm ESS}}}\frac{\alpha_{i}\zeta_{i,\tau}^{{\rm ESS}}}{0.8E_{i}^{\text{cap}}}
\end{array}\right)
\end{array}\right)
\end{align*}
}}%
which is subject to the constraints in (\ref{cons: C_ESS}). We have used
the revenue and cost expressions introduced in the
previous subsections and also the equalities in constraints (\ref{cons: service linkage a})
and (\ref{cons: service linkage b}) to deduce the TPN.

Since the decision variables $\{p_{i,\tau}^{{\rm br+fs,d}},p_{i,\tau}^{{\rm fs,c}}\}_{i\in\mathcal{N},\tau\in\mathcal{H}_{t}}$
do not appear in ${\rm TNP}(t)$, we can eliminate these redundant variables and simplify constraints (\ref{cons: service linkage a}) and
(\ref{cons: service linkage b}) into the following:
\begin{gather}
p_{i,t}^{{\rm re,c}}+p_{i,t}^{{\rm fr,c}}\le p_{i,t}^{{\rm c}}\le v_{i,t}^{{\rm c}}p_{i,\max}^{{\rm c}},\label{cons: service linkage a-1}\\
p_{i,t}^{{\rm fr,d}}\le p_{i,t}^{{\rm d}}\le(1-v_{i,t}^{{\rm c}})p_{i,\max}^{{\rm d}}-p_{i,t}^{{\rm sr,d}}+z_{i,t},\label{cons: service linkage b-1}
\end{gather}
for each $i\in\mathcal{N}$. Here the variable $z_{i,t}$ is an equivalent of $v_{i,t}^{{\rm c}}p_{i,t}^{{\rm sr,d}}$, satisfying constraint (\ref{cons: auxiliary z}). By minimizing
${\rm TNP}(t)$ with the new constraints, the solution of $(p_{i,\tau}^{{\rm br+fs,d}},p_{i,\tau}^{{\rm fs,c}})$
can then be recovered from (\ref{cons: service linkage a}) and (\ref{cons: service linkage b}).

Consequently, the ESS management
problem is defined as
\begin{align*}
\textbf{P0: } & \min-\text{TNP}(t)\\
\text{subject to, } & \eqref{cons: ES charge-discharge limits}-\eqref{cons: bill reduction},\thinspace\thinspace\eqref{cons: service linkage c}-\eqref{cons: service linkage b-1},\thinspace\thinspace\forall\thinspace i\in\mathcal{N},\tau\in\mathcal{H}_{t}
\end{align*}%
where the subscript $\tau$ instead of $t$ is used in all
constraints, and constraint (\ref{cons: C_ESS}) refers only to the
inequalities. In $\textbf{P0}$, the power variables are
$p_{i,\tau}^{{\rm c}}$, $p_{i,\tau}^{{\rm re,c}}$ and $p_{i,\tau}^{{\rm fr,c}}$
for charge and $p_{i,\tau}^{{\rm d}}$, $p_{i,\tau}^{{\rm fr,d}}$
and $p_{i,\tau}^{{\rm sr,d}}$ for discharge of each ESS $i\in\mathcal{N}$
in each time slot $\tau\in\mathcal{H}_{t}$, and $p_{\tau}^{{\rm re,sc}}$
and $p_{\tau}^{{\rm re,s}}$ for the customer to self-consume and
sell available renewable energy in each time slot $\tau\in\mathcal{H}_{t}$.
The auxiliary variables are binary variables $\{v_{i,\tau}^{{\rm c}}\}_{i\in\mathcal{N}}$,
$v_{\tau}^{{\rm fr}}$ and $v_{\tau}^{{\rm sr}}$
for all $\tau\in\mathcal{H}_{t}$, and real variables $\{\zeta_{i,\tau}^{{\rm ESS}}\}_{i\in\mathcal{N}}$
and $z_{i,\tau}$ for all $\tau\in\mathcal{H}_{t}$ .

The objective and constraints of $\textbf{P0}$ are linear in the variables, except constraint $\eqref{cons: C_ESS}$
which is quadratic in the decision variables $\{p_{i,\tau}^{{\rm c}},p_{i,\tau}^{{\rm d}}\}$ for each $i\in\mathcal{N}$ and $\tau\in\mathcal{H}_{t}$.
Therefore, the problem is a mixed-integer quadratic program
(MIQP) and can be solved by standard MIQP solvers. Once an optimal solution
is obtained, only the part for the first time slot $t$ will be
implemented to dispatch the ESSs. The schedule for the future
time slot $(t+1)$ will be determined in a similar way by shifting the
time horizon forward by one slot and then solving the new optimization
problem.

\section{Case Study \label{sec: Case-studies}}

This section assesses the economics of multi-use ESSs based on the scheduling approach developed above.

\subsection{Simulation setup and input data}

The  demand is scaled historical
hourly electricity demand of a college in California for a summer
week \cite{consumptiontrace2002}. The solar-PV generation, which has a peak equal to 60\% of the peak demand, is
scaled historical hourly generation power for a summer week in Brussels \cite{Solar-PV-ELIA2013}. The customer deploys two ESSs with specifications given in Table \ref{tb: ESS parameters}.
Their charge and discharge aging costs are estimated
by (\ref{cons: C_ESS}), in which the model parameters are set the same as those in \cite{hu2016towards}.

\begin{table}
\caption{\textsc{Parameters of the ESSs. The Power is in Unit of kW and the Energy is in Unit of kWh.}}
\label{tb: ESS parameters}

\centering{}%
\begin{tabular}{cccccccc}
\hline
Type & $E_{i}^{\text{cap}}$ & $s_{i,\min}$ & $s_{i,\max}$ & $p_{i,\max}^{{\rm c}}$ & $p_{i,\max}^{{\rm d}}$ & $\eta_{i}^{{\rm c}}$  & $\eta_{i}^{{\rm d}}$\tabularnewline
\hline
1 & 480 & 0.2 & 0.9 & 102 & 74 & 0.82 & 0.88\tabularnewline
2 & 720 & 0.2 & 0.9 & 148 & 113 & 0.85 & 0.90\tabularnewline
\hline
\end{tabular}
\vspace{-18pt}
\end{table}

The hourly regulation signal and associated market
clearing price are from the real operational records of PJM \cite{ancillary_src_data2015}. So are the hourly spinning
reserve prices. The price of purchasing electricity from the market is obtained
from PG \& E \cite{hu2016towards}, which consists of peak, mid-peak and off-peak prices for different periods of a day. The price
of selling electricity to the market is 60\%
of the purchase price. The purchase and sale powers are unrestricted in our study.

\subsection{Profitability of the multi-use ESSs}

\begin{figure}
\noindent \begin{centering}
\includegraphics[scale=0.395]{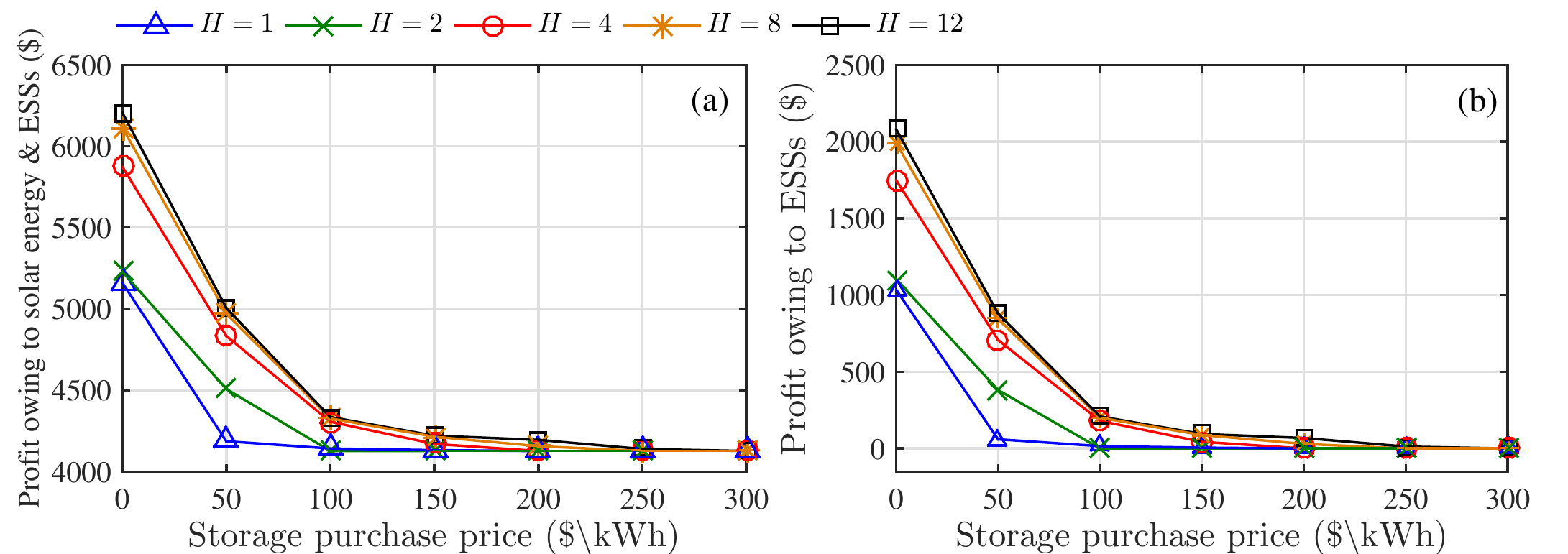}
\par\end{centering}

\caption{{\small{}Net profit vs. the rolling (look-ahead) horizon size.}}
\label{fig: profit-vs-horizon-size}
\vspace{-6pt}
\end{figure}

When the input data within the horizon is perfectly known, the total net profit obtained is shown
in Fig. \ref{fig: profit-vs-horizon-size}(a). The profit decreases
with the storage purchase price, and vanishes once the purchase price
is higher than 300 \$/kWh. On the other
hand, the profit increases with the horizon size $H$, but the marginal benefit decreases and becomes small once the horizon
size is larger than 4.

The above profit contains the contribution of solar energy which is
assumed free here. To get rid of this contribution and obtain the value solely for the use of ESSs, we subtract the above profit
with the one obtained without having any ESSs. This yields the reduced profits
shown in Fig. \ref{fig: profit-vs-horizon-size}(b).

\begin{figure}
\noindent \begin{centering}
\includegraphics[scale=0.495]{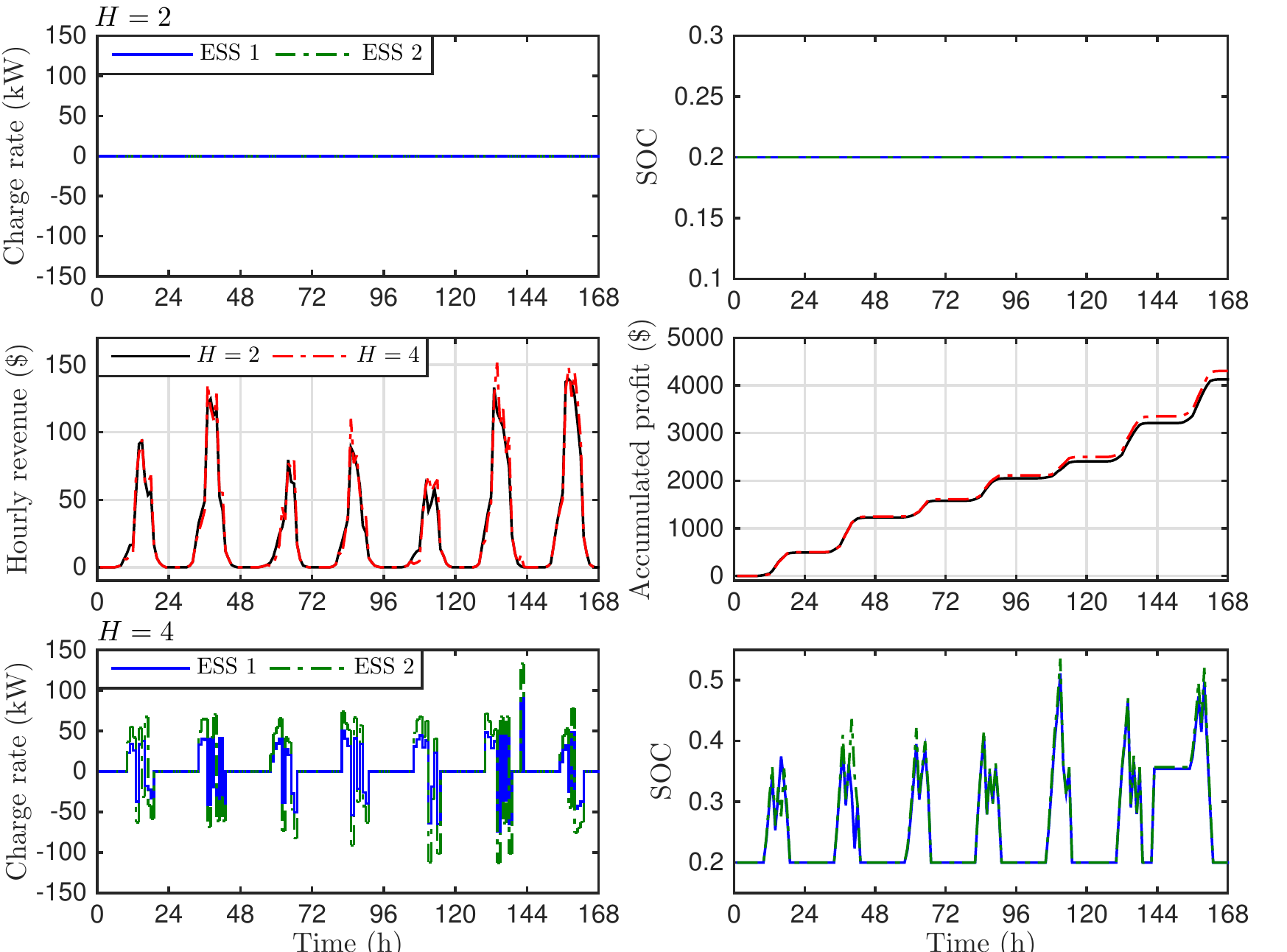}
\par\end{centering}

\protect\caption{{\small{}Battery dynamics and the resulting revenues and profits.
\label{fig: detailed-schedules}}}
\vspace{-16pt}
\end{figure}

When the storage purchase price is fixed to 100 \$/kWh, the charge/discharge schedules for $H$ equal to 2 and 4 are shown in Fig. \ref{fig: detailed-schedules}. The associated revenues and storage operating costs are given in Table \ref{tb: revenues-vs-costs}. As can be seen, the higher profit for the case with $H=4$ owes to appropriate use of the ESSs for regulation and reserve services and for reducing the electricity
bill. The results in Table \ref{tb: revenues-vs-costs} also indicate
that, in the absence of free solar energy, using ESSs to support
only one or two of the services may not cover the associated
cost and hence would be unable to yield a positive profit.

\begin{table}
\vspace{-16pt}
\caption{\textsc{Revenues and Costs in Unit of \$.}}
\label{tb: revenues-vs-costs}

\centering{}%
\begin{tabular}{ccccccc}
\hline
$H$ & $R_{\text{sc}}$ & $R_{\text{fr}}$ & $R_{\text{sr}}$ & $R_{\text{br}}$ & Storage aging cost & Net profit\tabularnewline
\hline
2 & 4126 & 0 & 0 & 0 & 0 & 4126\tabularnewline
4 & 4126 & 506 & 81 & 159 & 569 & 4303\tabularnewline
\hline
\end{tabular}
\vspace{-12pt}
\end{table}

\subsection{Impact of forecast errors}

\begin{figure}
\noindent \begin{centering}
\includegraphics[scale=0.39]{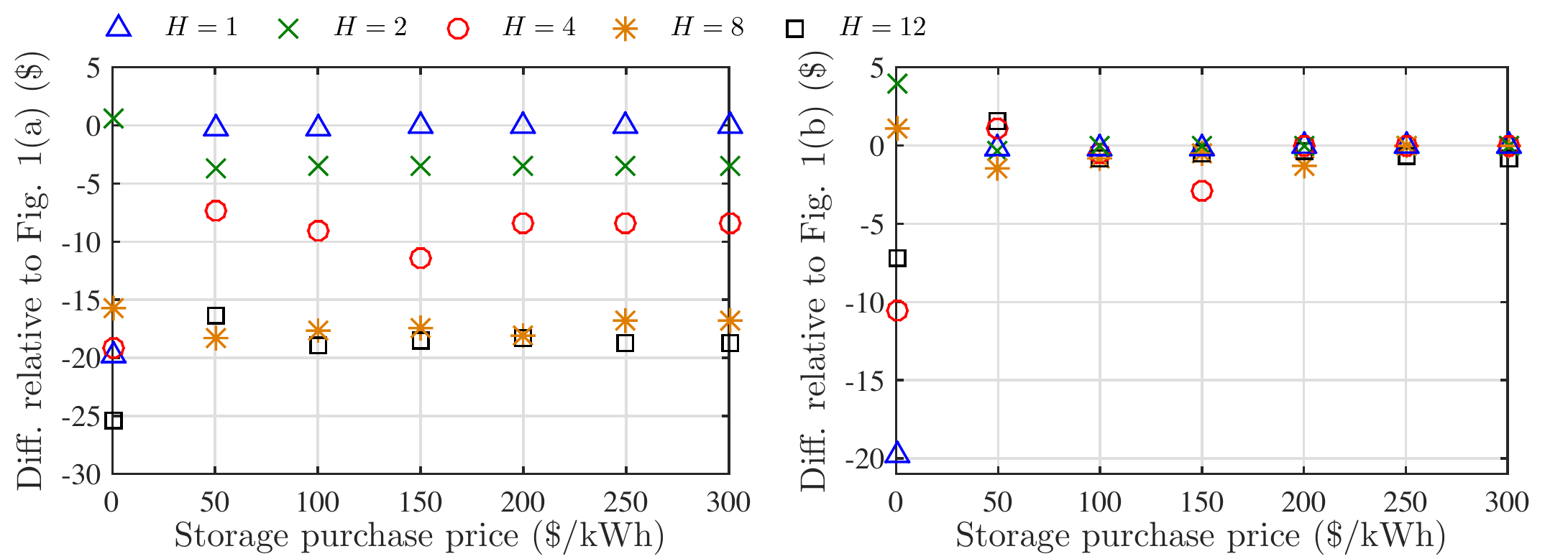}
\par\end{centering}

\protect\caption{{\small{}Profit differences relative to the profits in Fig. \ref{fig: profit-vs-horizon-size}. \label{fig: sensitivity-to-forecast-errors}}}
\vspace{-16pt}
\end{figure}

It is of interest to see how forecast errors affect the economic results. Let the load demand, the renewable generation power, the regulation market clearing price and the spinning reserve availability price be forecasted with zero-mean and uniformly distributed errors. The maximum errors are proportional to the magnitude changes in the true data for two sequential time slots, and the proportion coefficients increase with the forecast horizon. Further the forecast is capped within 80\% of the minimum and 120\% of the maximum true values. The differences of the resulting net profits relative to those in Fig. \ref{fig: profit-vs-horizon-size}(a)-(b) are shown in Fig. \ref{fig: sensitivity-to-forecast-errors}. As observed, the differences are mostly negative (it can be positive as rolling-horizon optimization may be sub-optimal in the long run), indicating losses of profits caused by the forecast errors. Nevertheless, the magnitudes are small relative to the reference profits. This indicates that the economic assessment approach is somehow robust to forecast errors.

\section{Conclusions}  \label{sec: Conclusions}

This paper developed a rolling-horizon optimization approach to schedule customer-sited ESSs for multi-service provision. The operating cost and yielded revenues were used
to assess the economics of the ESSs. The effectiveness of the proposed approach was illustrated with case studies. Future research will investigate the impact of storage energy and power capacities on the economics of multi-use ESSs.

\bibliographystyle{IEEEtran}
\bibliography{HuWangGooi_ESS_TENCON2016}
\end{document}